\documentclass{article}
\usepackage[utf8]{inputenc}
\usepackage{graphicx}
\usepackage{caption}
\usepackage{subcaption}
\usepackage{amsmath}
\usepackage{amsfonts}
\begin{document}
\title{A proof using B\"ohme's Lemma that no Petersen family graph has a flat embedding}
\author{Joel Foisy, Catherine Jacobs, Trinity Paquin, Morgan Schalizki, Henry Stringer}
\date{January 2023}

\maketitle
\begin{abstract}

Sachs and Conway-Gordon used linking number and a beautiful counting argument to prove that every graph in the Petersen family is intrinsically linked (have a pair of disjoint cycles that form a nonsplit link in every spatial embedding) and thus each family member has no flat spatial embedding (an embedding for which every cycle bounds a disk with interior disjoint from the graph). We give an alternate proof that every Petersen family graph has no flat embedding by applying B\"{o}hme's Lemma and the Jordan-Brouwer Separation Theorem. 
\end{abstract}

\section{Introduction}

%In this paper we will be exploring a new proof of how the Petersen Family Graphs such as $K_6$ and $P_7$ are not flat.

%\begin{figure*}[ht!]
   % \centering
  %  \begin{subfigure}[t]{0.5\textwidth}
     %   \centering
        %\includegraphics[width=.525\linewidth]{k_6.png}
        %\captionof{figure}{$K_6$}
        %\label{fig:test%1}
   % \end{subfigure}%
   % ~ 
   % \begin{subfigure}[t]{0.5\textwidth}
     %   \centering
        %\includegraphics[width=.4\linewidth]{p7.jpg}
        %\captionof{figure}{$P_7$}
        %\label{fig:test^2}
   %^ \end{subfigure}
%^\end{figure*}

Sachs \cite {sachs0}, \cite{sachs} and Conway-Gordon \cite{cg} (for $K_6$) proved, using the linking number and a beautiful counting argument, that every graph in the Petersen family is intrinsically linked (that is, has a pair of disjoint cycles that form a nonsplit link in every spatial embedding). As a consequence of being intrinsically linked, every graph in the Petersen family has no flat embedding. A 
\textit{flat} spatial embedding is one in which every cycle bounds a disk with interior disjoint from the graph. When such a disk is attached to a cycle, we say that cycle has been \textit{paneled}.

We have devised a new family of proofs that every Petersen family graph has no flat embedding, without using intrinsic linking. We sketch our methods here. Given a graph $G$ with a flat embedding, consider a set of cycles in $G$ such that, pair-wise, intersections are either connected or empty. We call such a set of cycles a \textit {B\"ohme system}. B\"ohme's Lemma \cite {bohme} concludes that for a B\"ohme system of cycles in a flat embedding of $G$, it is possible to panel the cycles in the system simultaneously--that is, each cycle bounds a disk with interior disjoint from the graph and from all other disks. The Jordan-Brouwer Separation Theorem then allows us to conclude that any sphere resulting from paneling a B\"ohme system (which we will call a \textit{B\"ohme sphere}) will separate space into an inside (bounded) and an outside (unbounded) component. Since we work in the PL category, by Alexander's result \cite{alex}, the inside component is a topological ball.

For the sake of contradiction, we suppose a Petersen family graph has a flat embedding and then find an appropriate B\"ohme system. We then argue that some resulting B\"ohme sphere must be intersected by an edge of the embedded graph, which contradicts flatness. In this way, we prove that each Petersen family graph cannot have a flat embedding. In the next section, we provide further details for each graph in the Petersen Family. 

%\section{Background}

%A \emph{spatial embedding} of a graph is a way to place a graph in space, so that vertices are points and edges are arcs that meet only at vertices. 

We conclude this introductory section with some terminology and background. We work in the PL category throughout. 

Suppose a graph $G$ has a three cycle with vertices $a,b,c$. To perform a $\Delta$ -Y exchange on the three-cycle in $G$, an extra vertex $y$ is added, edges $ab,bc,$ and $ca$ are removed, and the edges $ay,by,$ and $cy$ are added to transform $G$ into $G'$. If a graph $H$ has a degree 3 vertex, $y$, to perform a Y-$\Delta$ exchange on $H$ at $y$, we begin with the Y, having vertices $a,b,c, y$ and edges $ay,by,$ and $cy$. We remove the $y$ vertex and all edges incident to $y$ and then add edges $ab,bc,$ and $ca$ to make a triangle and transform $H$ to the graph $H'$. 

%\begin{center}
%\includegraphics[scale=.75]{Triangle-Y.png}
%\end{center}
 
The Petersen family of graphs consists of $K_6$ and $K_{3,3,1}$ and the graphs obtained from $K_6$ and $K_{3,3,1}$ by $\Delta-Y$ exchanges (equivalently, from $K_6$ by $\Delta-Y$ and $Y-\Delta$ exchanges). We will denote the remaining graphs in this family, with subscripts indicating the number of vertices, as $P_7$, $P_8$, $P_9$, $P_{10}$ (which is the classic Petersen graph) and $K_{4,4}-e$. Robertson, Seymour and Thomas proved Sachs' conjecture that the Petersen family of graphs form the complete minor-minimal set of graphs that are intrinsically linked. They further showed that a graph has a flat embedding if and only if it  has a linkless embedding \cite{rst}.

%\item {\bf Spherically flat graphs.} First, a spatial graph $f(G)$ is \emph{flat} if it is possible to attach a topological disk to every cycle, $C$, of $G$ so that the disk intersects $C$ along its boundary, and the interior of the disk is disjoint from $f(G)$. Every flat embedding is linkless. Robertson, Seymour and Thomas showed that a graph has a flat embedding if and only if it  has a linkless embedding. We'll call a graph $G$ \emph{flat} if it has a flat spatial embedding (Robertson, Seymour, and Thomas).\\

%A graph $G$ is \emph{spherically flat} if there is a flat spatial embedding of $G$, $f(G)$, such that every planar subgraph $P$ of $G$ lies on an embedded sphere $S$ that intersects $f(G)$ only in $f(P)$. By the Jordan-Brouwer theorem, any vertex or edge in $f(G)-f(P)$ lies on one side of $S$ (for an edge, its endpoints lies in $S$).\\

%Note that if $G$ is spherically flat, then $G$ is flat (every cycle is itself a planar subgraph, so lies on a sphere, and from there, take a bounding disk) and $G$ has a linkless embedding.\\

\section{Proofs}

In this section, we go through details of the proof for each graph in the Petersen family.
%\begin{itemize}

 % \centering
   % \begin{subfigure}[t]{0.5\textwidth}
      %  \centering
       % \includegraphics[width=.525\linewidth]{k_6.png}
       % \captionof{figure}{$K_6$}
       % \label{fig:test1}
   % \end{subfigure}%
   % ~ 
   % \begin{subfigure}[t]{0.5\textwidth}
      %  \centering

\subsection {$K_6$}

\begin{figure}[ht!]
   
       \hskip 1.7in \includegraphics[width=.23\linewidth]{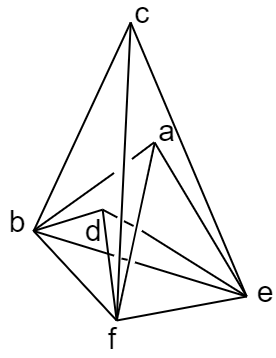}
        \caption{A spatial embedding of the subgraph of $K_6$ formed by $\mathcal{C}$.}
  
   % \end{subfigure}
 
    \label{k6}
\end{figure}

%\ref{test}
Denote the vertices of $K_6$ as $\{a,b,c,d,e,f\}$.
By way of contradiction suppose that $G = K_6$ has a flat embedding. Consider the collections of cycles:
\begin{center}
    $C_1 = \{(b,e,f), (b,c,e), (c,e,f), (b,c,f)\}$\\
    $C_2 = \{(b,e,f), (a,b,e), (a,e,f), (a,b,f)\}$\\
    $C_3 = \{(b,e,f), (b,d,e), (d,e,f), (b,d,f)\}$\\
\end{center}
    Let $\mathcal{C} = \cup^3_{i=1}C_i$. Note that the cycles above are induced cycles where $C_1 \subset G[\{b,c,e,f\}], C_2 \subset G[\{a,b,e,f\}]$, and $C_3 \subset G[\{b,d,e,f\}]$, where $G[\{a,b,c,d\}]$ is the subgraph of $K_6$ induced by the vertices $a,b,c,d$. Notice  $(b, e, f) \in C_i$ for $i = 1, 2, 3$, see Figure \ref{k6}.
    
    %$D_1, ..., D_{10}$. Let $T_i = C_i \cup (\cup^{k}_j D_j)$ where $i = 1, 2, 3$, $j = 4i - 3$ and $k = j + 4$. 
    
   Take a flat embedding of $K_6$. As all the cycles of $\mathcal{C}$ have connected pairwise intersections, we can apply B\"{o}hme's Lemma to get simultaneous paneling disks of the ten 3-cycles that make up $\mathcal{C}$. As we panel each $C_i$, we form a sphere $T_i$, then by Jordan-Brouwer Theorem applied to, say, $T_1$, there are inside and outside components, $B_1$ and $B_2$ respectively, with $B_1 \cup B_2={\mathbb R}^3-T_1$, and, say, $B_1$ a ball (which follows from Alexander's result \cite{alex} in the PL category). Without loss of generality as the edge $\overline{ad}$ is disjoint from $T_1$, we may assume that $ad \subset B_1$. As $ad \subset B_1$ and $a \in T_2$, $d \in T_3$ then $T_2, T_3 \subset B_1 \cup T_1$. This means that $T_2$ and $T_3$ are in the closure of $B_1$ (this is what we mean when we later say that $T_2$ and $T_3$ lie inside of $T_1$). Now similarly observe there exists, via Jordan-Brouwer, an inside and outside region for $T_2$. Without loss of generality suppose that $T_3$ is in the closure of the inner region (ball) of $T_2$. As $d \not\in T_2$ and $d \in T_3$, $d$ is in the inner ball of $T_2$. As $T_2$ is inside $T_1$ and $c \in T_1$, c is outside $T_2$. As edge $cd$ exists and as c and d are on opposite sides of $T_2$, edge $cd$ must intersect a disk of $T_2$.  This brings us a contradiction. Thus $K_6$ has no flat embedding.

\noindent
{\bf The remaining cases}\\
For all remaining Petersen family graphs except $P_{10}$ and $K_{4,4}-e$, similarly to the proof for $K_6$, we consider a flat embedding, and a cycle $C$ with vertex set $V(C)$, together with three vertices disjoint from $V(C)$; $v_1$, $v_2$ and $v_3$, such that each $G[V(C) \cup \{v_i\}]$ forms a B\"{o}hme system. After panelling the B\"{o}hme system, the result is three spheres that intersect only along $C$ and the disk that panels $C$. In each case, two vertices are on opposite sides of a resulting B\"{o}hme sphere. In each case, the two vertices are connected either by an edge, or through a $Y$ that is disjoint from $G[V(C) \cup \{v_i\}]$ for the appropriate $i$, and the result is a contradiction of flatness. More details follow.

%\item $P_7$
\subsection {$P_7$}
Denote the six vertices of $K_6$ by $\{a,b,c,1,2,3\}$, and suppose that a $\Delta -Y$ exchange took place between $1,2,$ and $3$. By way of contradiction, suppose that $P_7$ has a flat embedding. Consider the induced subgraphs:\begin{center}
   $G_1 = G[\{a,b,c,1\}]$\\
     $G_2 = G[\{a,b,c,2\}]$\\
  $G_3 = G[\{a,b,c,3\}]$
\end{center}

%The sets of induced cycles in each induced subgraphs are:\\ $\mathcal{B}_1 = \{(1,a,2,v),(1,a,2,b),(1,v,2,b)\}$, $\mathcal{B}_2 = \{(1,a,2,v),(1,a,2,c),(1,v,2,c)\},$ and $\mathcal{B}_3= \{(1,a,2,v),(1,a,3,v),(2,a,3,v)\}$.
Note that each $G_i$ contains the cycle $(a,b,c)$. One can check that the induced cycles in $G_1 \cup G_2 \cup G_3=K_{3,1,1,1}$ form a B\"ohme system. Take a flat embedding of $P_7$. The simultaneous paneling of the B\"ohme system turns each $G_i$ into a B\"ohme sphere. From here, the argument follows analogously to the argument for $K_6$, except the vertices $\{1,2,3\}$ are connected by a $Y$ instead of a triangle.
\begin{figure}[h]
\begin{center}
\includegraphics[scale=.48]{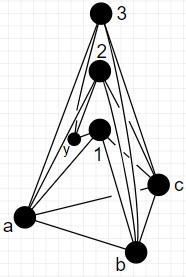}
\end{center}
\caption{A spatial embedding of $P_7$.}
\label{P7'}
\end{figure}

%Let $S_i$ denote the sphere formed by $G_i$. As $3 \notin G_1$, without loss of generality, we can assume that $3$ is outside of $S_1$.
%As $2$ and $3$ are connected through $Y$, and neither 2, nor 3, nor $Y$ lie in $S_1$, then $2$ and $3$ must also both be outside $S_1$. By similar reasoning, either both $1$ and $2$ are inside $S_3$, or both $1$ and $3$ are inside $S_2$. Without loss of generality, we will suppose both $1$ and $2$  are inside $S_3$ (as suggested in Figure \ref{P7'}). It follows that $1$ is inside $S_2$ and $3$ is outside of $S_2$. Since $1$ is connected to $3$ through $Y$, it follows that the path 1Y3 must intersect a panel of $S_2$. This contradicts $P_7$ having a flat embedding.

%Old
%As $b$ and $3$ are connected by an edge and $3$ is on the outside of $S_2$, $b$ is also on the outside of $S_2$. But as $3$ is on the inside of the sphere formed by $\mathcal{B}_2$, $b$ is on the outside of the sphere formed by $\mathcal{B}_2$, which is a contradiction. See Figure \ref{P7'}. Thus, $P_7$ is not flat.

%\item $K_{4,4} - e$ 
\subsection{$K_{4,4} - e$}
Denote the vertices of $K_{4,4}-e$ as $V_1=\{a,b,c,d\}$ and $V_2\{1,2,3,4\}$ where all possible edges from $V_1$ to $V_2$ are present, except $e= 4a$\\

%Looking back at the embedding of $P_7$ from Figure \ref{P7}, we can get to $K_{4,4}-e$ by exchanging the triangle that is disjoint from the $Y$. If we start the B\"ohme system described for $P_7$, we get a similar system for $K_{4,4}-e$, and a similar contradiction.

%Here are more details. 
\begin{figure}[h]
\begin{center}
    \includegraphics[scale=.7]{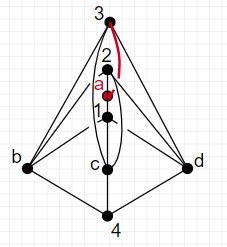}
\end{center}
\caption{A spatial embedding of $K_{4,4}-e$. Possibly the edge $a3$ intersects a panel.}
\label{k44}
\end{figure}

Assume by way of contradiction that $K_{4,4} - e$ has a flat embedding. 

%Then by definition of flat, there is a spacial embedding, such that for each cycle there is a disk bounded by the cycle with an interior disjoint from the rest of the graph. 

Consider the set of cycles $A_i$ listed bellow, where $i \in \{1,2,3\}$.  \begin{center}
   $A_1 = \{(1,b,4,c),(1,c,4,d),(1,b,4,d)\} \subseteq G[1,4,b,c,d]$\\
 $A_2 = \{(2,b,4,c),(2,c,4,d),(2,b,4,d)\} \subseteq G[2,4,b,c,d]$\\
     $A_3 = \{(3,b,4,c),(3,c,4,d),(3,b,4,d)\} \subseteq G[3,4,b,c,d]$

\end{center} 
 Notice that $A_1 \cup A_2 \cup A_3$ satisfies the hypotheses of B\"ohme's Lemma with edges forming $K_{4,3}$. By B\"ohme's lemma, we can simultaneously panel each of these cycles. Note that all paneled cycles pass through the vertex $4$. Denote the sphere resulting from paneling $A_i$ as $S_i$, where $i \in \{1,2,3\}$. Note that each $A_i$ contains the $Y$ on vertices $4, b, c, d$, as well as the vertex $i$. From here, the proof follows similarly to the argument used for $P_7$ ($K_6$), with the base cycle (triangle) $C$ replaced with the $Y$, where for $i\neq j$, $S_i \cap S_j$ is precisely the $Y$.

 %Note that since 1, 2 and 3 are all connected to $a$ by an edge, the following is without loss of generality. We may take $1$ and $2$ both inside $S_3$. Since 1 and 2 are  connected to $a$ by an edge, by the Jordan-Brouwer Theorem $a$ must exist outside $S_3$, and moreover $a$ must exist between $S_1$ and $S_2$, or else edge $\overline{a1}$ intersects $S_2$ or edge $\overline{a2}$ must intersect $S_1$. However, 3 is also connected to $a$ by an edge, and $3$ lies outside of both $S_1$ and $S_2$ so the edge between 3 and $a$ must intersect $S_1$ or $S_2$, which contradicts flatness.  Therefore $K_{4,4} - e$ has no flat embedding. \\
 
 %%Is there another inequivalent B\"ohme system for $K_{4,4}-e$? There would not appear to be so, but this should be checked.

%\item $K_{3,3,1}$
\subsection{$K_{3,3,1}$}
%\vskip -.7in
\begin{figure}[h]
%\begin{center}
   % \includegraphics[scale=.73]{K_3,3,1.png}
 %\hskip 1.65in  % \includegraphics[scale=.36]{k331.pdf}
 
\hskip 1.5in \includegraphics[scale=.5]{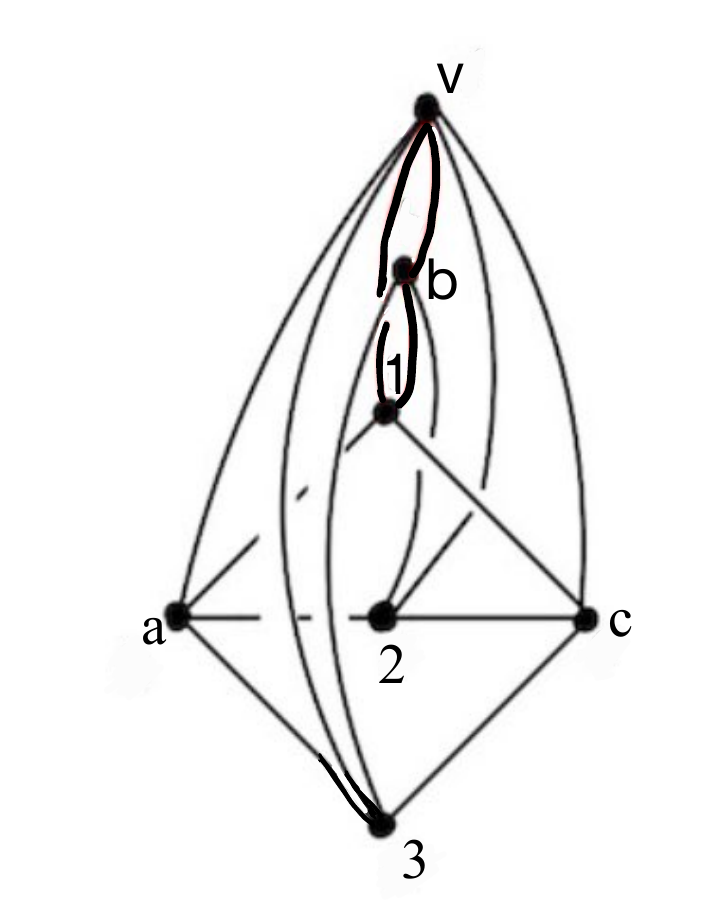}
%\end{center}
%\vskip -.7in
\caption{A spatial embedding of $K_{3,3,1}$.}
\label{k331}
\end{figure}

Label the vertices of $K_{3,3,1}$ as in Figure \ref{k331}.
Assume by way of contradiction that $K_{3,3,1}$ has a flat embedding. Consider the collection of cycles:
\begin{center}
   $\mathcal{S}_1 = \{(a,2,c,3),(a,1,c,2),(a,1,c,3)\} \subseteq G[a,2,c,3,1]$\\
   $\mathcal{S}_2 = \{(a,2,c,3),(3,b,2,c),(3,b,2,a)\} \subseteq G[a,2,c,3,b]$\\ 
    $\mathcal{S}_3 = \{(a,2,c,3),(c,2,v),(c,v,3),(3,v,a), (a,v,2)\} \subseteq G[a,2,c,3,v]$
\end{center}

Notice that $\mathcal{S}_1 \cup \mathcal{S}_2 \cup \mathcal{S}_3$ satisfies B\"ohme's lemma, and after paneling each, $\mathcal{S}_i$ forms a B\"ohme sphere. These B\"ohme spheres share precisely the same base cycle $(a,2,c,3)$ with the corresponding panelling disk. From here, the argument follows similarly to the argument for $K_6$.

%Without loss of generality assume the vertex $b$ lies outside of the sphere generated by $\mathcal{S}_1$. Then since b is outside of the sphere generated by $\mathcal{S}_1$, the sphere generated by $\mathcal{S}_2$ must lie in the closure of the region outside of the sphere generated by $\mathcal{S}_1$. Furthermore, since v and b are connected by an edge, by the Jordan-Brouwer Theorem, they must both be inside or both be outside of the sphere generated by $\mathcal{S}_1$. Without loss of generality assume they are both outside of the sphere generated by $\mathcal{S}_1$. Then since $v$ is outside of the sphere generated by $\mathcal{S}_1$, the sphere generated by $\mathcal{S}_3$ must lie in the closure of the region outside of the sphere generated by $\mathcal{S}_1$. Since $b$ and $1$ are connected by an edge, they must both lie inside the sphere determined by $\mathcal{S}_3$. However, consider the edge between $v$ and $1$. The vertex $v$ is outside of the sphere generated by $\mathcal{S}_2$, however, $1$ is inside the sphere generated by $\mathcal{S}_2$. However, since they are connected by an edge they both must lie inside or outside the sphere generated by $\mathcal{S}_2$ by the Jordan-Brouwer theorem. Thus we have a contradiction and $K_{3,3,1}$ has no flat embedding.

%\item $P_8$\\
\subsection{$P_8$}
%Looking back at the B\"ohme system above for $K_{3,3,1}$, there are, up to equivalences, three different triangles that can be exchanged to yield $P_8$:

%Similarly, looking at $P_7$, ...

Denote the vertices of $P_8$ by $\{v, a,b,c, 1,2,3, y\}$, where $\{v\}$, $\{a,b,c\}$, and $\{1,2,3\}$ are the partitions in $K_{3,3,1}$, and triangle $(v,a,1)$ was exchanged to a $Y$, with resulting new $Y$ vertex $y$ in the resulting $P_8$.

\begin{figure}[h]
  %  \centering
   % \includegraphics{}
   % \caption{Caption}
   % \label{fig:my_label}
%\end{figure}
\begin{center}
\includegraphics[scale=.7]{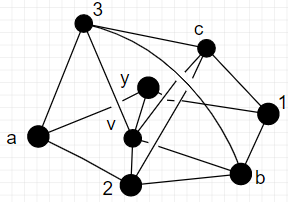}
\end{center}
\caption{A spatial embedding of $P_8$.}
\label{P8}
\end{figure}
 By way of contradiction, suppose that $P_8$ is has a flat embedding. Consider the induced subgraphs:
\begin{center}
    \item $G_1 = G[a,y,1,b,2,v]$\\
    \item $G_2 = G[a,y,1,b,2,c]$\\
    \item $G_3 = G[a,y,1,b,2,3]$
\end{center}
and their respective B\"ohme spheres formed by panelling their induced cycles (which is possible by B\"{o}hme's Lemma):
\begin{center}
    $\mathcal{B}_1 = \{(a,y,1,b,2),(v,2,a,y),(v,y,1,b),(v,b,2)\}$\\
     $\mathcal{B}_2 = \{(a,y,1,b,2),(a,y,1,c,2),(b,1,c,2)\}$\\
     $\mathcal{B}_3 = \{(a,y,1,b,2),(a,y,1,b,3),(a,2,b,3)\}$
\end{center}
 Note that each $\mathcal{B}_i$ contains the cycle $(a,y,1,b,2)$. From here, the argument follows analogously to the argument for $K_6$.
 
 %Without loss of generality, as $v$ and $3$ are connected by an edge, by the Jordan-Brouwer Theorem, both $v$ and $3$ are on the outside of the sphere formed by $\mathcal{B}_2$. Then $v$ is on one side of the sphere formed by $\mathcal{B}_3$. Without loss of generality, suppose $v$ is inside of $\mathcal{B}_3$. Thus $v$ is inside $\mathcal{B}_3$ and outside $\mathcal{B}_2$. Thus every edge incident to $v$ has interior $\mathcal{B}_3$ and outside $\mathcal{B}_2$. It follows that $\mathcal{B}_1$ lies in the closure of the interior of the inside component of $\mathcal{B}_3$ and in the closure of the exterior of the outside component of $\mathcal{B}_2$. Finally, it follows that $c$ is inside $B_1$ and $3$ is outside $B_1$. However, the edge connecting $3$ and $c$ would then intersect $\mathcal{B}_1$, which contradicts flatness of the embedding. Thus, $P_8$ has no flat embedding.
 
 %%Then, as $v$ is connected to $c$ and $v$ is inside of $\mathcal{B}_3$, $c$ is inside of $\mathcal{B}_3$ by the Jordan-Brouwer Theorem. As $c$ is connected to $3$, 3 is inside of $\mathcal{B}_1$ by the Jordan-Brouwer Theorem. However, the edge connecting 3 and $v$ is both inside and outside of $\mathcal{B}_2$, which contradicts the Jordan-Brouwer Theorem. Thus, $P_8$ is not flat.

%\item $P_9$\\
\subsection{$P_9$}

\begin{figure}[h]
\begin{center}
    \includegraphics[scale=.72]{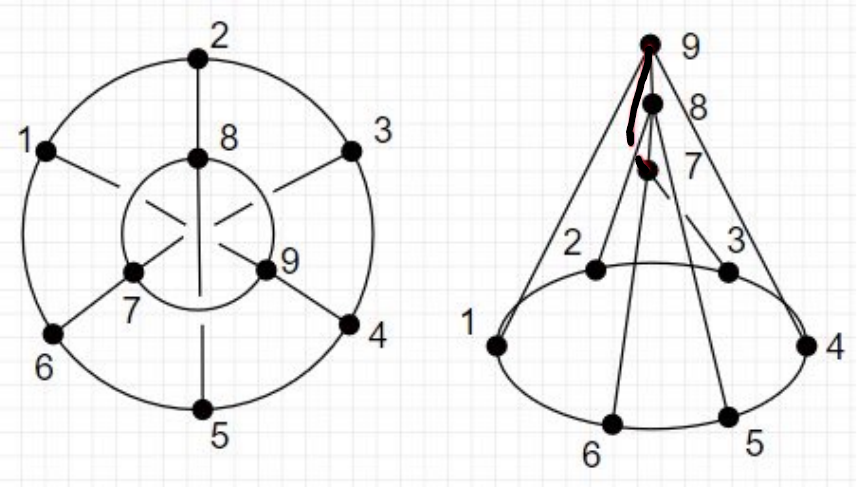}
    %p9
\end{center}
\caption{Spatial embeddings of $P_9$.}
\label{P9}
\end{figure}

Let $P_9$ be labelled as in Figure \ref{P9}.

Assume by way of contradiction that $P_9$ is has a flat embedding. Consider a flat embedding, and consider the sets of cycles ($S_i$ where $i \in \{1,2,3\}$): \begin{center}
    $S_1 = \{(1,2,3,4,5,6),(6,7,3,4,5),(6,7,3,2,1)\}\subset G[1,2,3,4,5,6,7]$ \\
    $S_2 = \{(1,2,3,4,5,6),(2,8,5,4,3),(2,8,5,6,1)\} \subset G[1,2,3,4,5,6,8]$ \\
    $S_3 = \{(1,2,3,4,5,6),(1,9,4,3,2),(1,9,4,5,6)\}\subset G[1,2,3,4,5,6,9]$
\end{center} 
 Note that each $S_i$ contains the cycle $C=(1,2,3,4,5,6)$ and notice that by B\"ohme's lemma each of these cycles can be simultaneously panelled so each $S_i$ determines a sphere for $i=1,2,3$ with, for $i \neq j$, the sphere determined by $S_i$ intersecting the sphere determined by $S_j$ only along $C$ and its corresponding panelling disk. From here, the proof follows analogously to the proof for $K_6$.

\subsection{$P_{10}$}

\vskip -.25in
\begin{figure}[h]
\begin{center}
    \includegraphics[scale=.9]{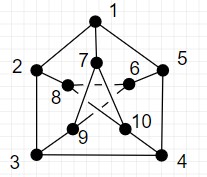}
\end{center}
\vskip -.1in
\caption{A spatial embedding of the Petersen graph.}
\label{Pete}
\end{figure}
Suppose by way of contradiction that $P_{10}$ has a flat embedding. Consider the induced subgraphs formed by each set of vertices below, where labels are from Figure \ref{Pete}.  \begin{center}
    $G_1 = G[\{1,2,3,4,5,6,9\}]$ \\
    $G_2 = G[\{1,2,3,4,5,7,10\}]$ \\
     $G_3 = G[\{1,2,3,4,5,6,8\}]$
\end{center} 
Consider a flat embedding of $P_{10}$. Then by B\"ohme's lemma, each 5-cycle and 6-cycle in $G_1,G_2,$ and $G_3$ can simultaneously bound disks to form three spheres, corresponding to each $G_i$. Note that the cycle $(1,2,3,4,5)$ is in each $G_i$.

    By the Jordan-Brouwer Theorem, since vertex 7 is connected to vertex 10 by an edge, they must both be outside or both be inside of the sphere formed by $G_1$. Without loss of generality, assume both vertices are outside. Then the sphere created by $G_2$ is in the closure of the region outside of the sphere created by $G_1$. Since vertex 8 is connected to both 6 and 10 by an edge, vertex 8 must be between the sphere created by $G_2$ and the sphere created by $G_1$. Thus the sphere formed by $G_3$ is sandwiched between the sphere created by $G_2$ and the sphere created by $G_1$ (that is, it lies in the closure of the outer region formed by $G_2$ and in the closure of the inner region formed by $G_1$). But now since 7 is connected to 9 by an edge, and 7 is in $G_2$ and 9 is in $G_1$, the sphere formed by $G_3$ intersects the edge connecting 7 to 9. Hence a paneled disk must have been punctured by edge 79, contradicting flatness. In Figure \ref{Pete}, the panelled cycle $(2,3,4,10,8)$ from $G_3$ would be punctured. Therefore $P_{10}$ has no flat embedding.
%\end{itemize}

We note that the $P_{10}$ case was similar to the earlier cases, except we needed a shared cycle $C$ and 5 additional vertices, and $G_1$ and $G_3$ shared the vertex $6$, in addition to the cycle $C$. The spheres they formed shared the base cycle (and disk) as well as the edge $5,6$. Instead of having a triangle or $Y$ connecting the extra vertices, here we have edges $8,10$ and $7,9$. Since $6$ lies in both $G_1$ and $G_3$ but not in $G_2$; $G_1$ and $G_3$ must be on the same side of $G_2$. Even with a swap of positions of $G_1$ and $G_3$ from our original proof, the edge $8,10$ would have intersected the B\"{o}hme sphere formed by $G_1$.

We could have done a $\Delta-Y$ exchange on $P_9$, on triangle $(7,8,9)$, to get $P_{10}$ with the same B\"ohme system we used on $P_9$ and the proof details analogous to those used on $P_9$.

 \section{Acknowledgements} This work was completed as part of a United States National Security Agency grant H89230-22-1-0008 Clarkson University: Summer (2022) Research Experience for Undergraduates in Mathematics.

\end{document}